\input amstex
\documentstyle{amsppt}
%
\catcode`@=11
\redefine\output@{%
  \def\break{\penalty-\@M}\let\par\endgraf
  \ifodd\pageno\global\hoffset=105pt\else\global\hoffset=8pt\fi  
  \shipout\vbox{%
    \ifplain@
      \let\makeheadline\relax \let\makefootline\relax
    \else
      \iffirstpage@ \global\firstpage@false
        \let\rightheadline\frheadline
        \let\leftheadline\flheadline
      \else
        \ifrunheads@ 
        \else \let\makeheadline\relax
        \fi
      \fi
    \fi
    \makeheadline \pagebody \makefootline}%
  \advancepageno \ifnum\outputpenalty>-\@MM\else\dosupereject\fi
}
\catcode`\@=\active
\nopagenumbers
\def\negskp{\hskip -2pt}
\def\Alpha{\operatorname{A}}
\def\divr{\operatorname{div}}
\def\vtrule{\vrule height 12pt depth 6pt}
\accentedsymbol\bbd{\kern 2pt\bar{\kern -2pt\bold d}}
\def\blue#1{#1}
\catcode`#=11\def\diez{#}\catcode`#=6
\catcode`_=11\def\podcherkivanie{_}\catcode`_=8
\def\mycite#1{\cite{\blue{#1}}\immediate\special{ps:
     ShrHPSdict begin /ShrBORDERthickness 0 def}}

\def\mytag#1{%
    \tag#1}
\def\mythetag#1{\thetag{\blue{#1}}\immediate\special{ps:
     ShrHPSdict begin /ShrBORDERthickness 0 def}}
\def\myrefno#1{\no#1}
\def\myhref#1#2{\blue{#2}\immediate\special{ps:
     ShrHPSdict begin /ShrBORDERthickness 0 def}}
\def\myEarXivlink{\myhref{http://arXiv.org}{http:/\negskp/arXiv.org}}
\def\myGeoCities{\myhref{http://www.geocities.com}{GeoCities}}
\def\mytheorem#1{\csname proclaim\endcsname{Theorem #1}}
\def\mythetheorem#1{\blue{#1}\immediate\special{ps:
     ShrHPSdict begin /ShrBORDERthickness 0 def}}
\def\mylemma#1{\csname proclaim\endcsname{Lemma #1}}

\def\mycorollary#1{\csname proclaim\endcsname{Corollary #1}}

\def\mydefinition#1{\definition{Definition #1}}

\font\eightcyr=wncyr8
\pagewidth{360pt}
\pageheight{606pt}
\topmatter
\title
On the Dirac equation in a gravitation field\\
and the secondary quantization.
\endtitle
\author
R.~A.~Sharipov
\endauthor
\address 5 Rabochaya street, 450003 Ufa, Russia\newline
\vphantom{a}\kern 12pt Cell Phone: +7-(917)-476-93-48
\endaddress
\email \vtop to 30pt{\hsize=280pt\noindent
\myhref{mailto:R\podcherkivanie Sharipov\@ic.bashedu.ru}
{R\_\hskip 1pt Sharipov\@ic.bashedu.ru}\newline
\myhref{mailto:r-sharipov\@mail.ru}
{r-sharipov\@mail.ru}\newline
\myhref{mailto:ra\podcherkivanie sharipov\@lycos.com}{ra\_\hskip 1pt
sharipov\@lycos.com}\vss}
\endemail
\urladdr
\vtop to 20pt{\hsize=280pt\noindent
\myhref{http://www.geocities.com/r-sharipov}
{http:/\negskp/www.geocities.com/r-sharipov}\newline
\myhref{http://www.freetextbooks.boom.ru/index.html}
{http:/\negskp/www.freetextbooks.boom.ru/index.html}\vss}
\endurladdr
\abstract
    The Dirac equation for massive free electrically neutral spin 1/2
particles in a gravitation field is considered. The secondary quantization
procedure is applied to it and the Hilbert space of multiparticle quantum 
states is constructed.
\endabstract
\endtopmatter
\loadbold
\loadeufb
\TagsOnRight
\document
\accentedsymbol\tbvartheta{\tilde{\overline{\boldsymbol\vartheta}}
\vphantom{\boldsymbol\vartheta}}

\rightheadtext{On the Dirac equation in a gravitation field \dots}
\head
1. The Dirac equation and its current. 
\endhead
    Let $M$ be a {\it space-time} manifold. It is a four-dimensional
orientable manifold equipped with a pseudo-Euclidean Minkowski-type 
metric $\bold g$ and with a {\it polarization}. The polarization of
$M$ is responsible for distinguishing the {\it Future light cone\/} 
from the {\it Past light cone\/} at each point $p\in M$ (see \mycite{1} 
for more details). Let's denote by $DM$ the bundle of Dirac spinors
over $M$ (see \mycite{2} and \mycite{3} for detailed description of
this bundle). In addition to the metric tensor $\bold g$ inherited from
$M$, the Dirac bundle $DM$ is equipped with four other basic spin-tensorial 
fields:
$$
\vcenter{\hsize 10cm
\offinterlineskip\settabs\+\indent
\vtrule
\hskip 1.2cm &\vtrule 
\hskip 5.2cm &\vtrule 
\hskip 2.8cm &\vtrule 
\cr\hrule 
\+\vtrule
\hfill\,Symbol\hfill&\vtrule
\hfill Name\hfill &\vtrule
\hfill Spin-tensorial\hfill &\vtrule\cr
\vskip -0.2cm
\+\vtrule
\hfill &\vtrule
\hfill \hfill&\vtrule
\hfill type\hfill&\vtrule\cr\hrule
\+\vtrule
\hfill $\bold g$\hfill&\vtrule
\hfill Metric tensor\hfill&\vtrule
\hfill $(0,0|0,0|0,2)$\hfill&\vtrule\cr\hrule
\+\vtrule
\hfill $\bold d$\hfill&\vtrule
\hfill Skew-symmetric metric tensor\hfill&\vtrule
\hfill $(0,2|0,0|0,0)$\hfill&\vtrule\cr\hrule
\+\vtrule
\hfill$\bold H$\hfill&\vtrule
\hfill Chirality operator\hfill&\vtrule
\hfill $(1,1|0,0|0,0)$\hfill&\vtrule\cr\hrule
\+\vtrule
\hfill$\bold D$\hfill&\vtrule
\hfill Dirac form\hfill&\vtrule
\hfill $(0,1|0,1|0,0)$\hfill&\vtrule\cr\hrule
\+\vtrule
\hfill$\boldsymbol\gamma$\hfill&\vtrule
\hfill Dirac $\gamma$-field\hfill&\vtrule
\hfill $(1,1|0,0|1,0)$\hfill&\vtrule\cr\hrule
}\quad
\mytag{1.1}
$$
As we see in the table \mythetag{1.1}, the metric tensor $\bold g$
is interpreted as a spin-tensorial field of the type $(0,0|0,0|0,2)$.
The Dirac bundle is a complex bundle over a real manifold. For this
reason spin-tensorial bundles produced from $DM$ are equipped with
the involution of complex conjugation $\tau$:
$$
\hskip -2em
\CD
@>\tau>>\\
\vspace{-4ex}
D^\alpha_\beta\bar D^\nu_\gamma T^m_n M@.
D^\nu_\gamma\bar D^\alpha_\beta T^m_n M.\\
\vspace{-4.2ex}
@<<\tau< 
\endCD
\mytag{1.2}
$$
Note that two fields $\bold g$ and $\bold D$ in \mythetag{1.1} are
real fields:
$$
\pagebreak
\xalignat 2
&\tau(\bold g)=\bold g,
&&\tau(\bold D)=\bold D.
\endxalignat
$$
Other fields $\bold d$, $\bold H$, and $\boldsymbol\gamma$ in the table
\mythetag{1.1} are not real spin-tensorial fields.
\mydefinition{1.1} A metric connection $(\Gamma,\Alpha,\bar{\Alpha})$ 
in $DM$ is a spinor connection real in the sense of the involution 
\mythetag{1.2} and concordant with $\bold d$ and $\boldsymbol\gamma$, i\.\,e\.
$$
\xalignat 4
&\nabla\bold d=0,
&&\nabla\boldsymbol\gamma=0,
&&\text{and}
&&\tau(\nabla\bold X)=\nabla(\tau(\bold X)),
\ \qquad
\mytag{1.3}
\endxalignat
$$
where $\bold X$ is an arbitrary smooth spin-tensorial field of the Dirac
bundle.
\enddefinition
\mytheorem{1.1} Any metric connection $(\Gamma,\Alpha,\bar{\Alpha})$
is concordant with all of the basic spin tensorial fields $\bold g$, 
$\bold d$, $\bold H$, $\bold D$, and $\boldsymbol\gamma$ listed in the
table \mythetag{1.1}.
\endproclaim
\noindent The theorem~\mythetheorem{1.1} means that from \mythetag{1.3} it
follows that 
$$
\xalignat 3
&\nabla\bold g=0,
&&\nabla\bold H=0,
&&\nabla\bold D=0.
\ \qquad
\mytag{1.4}
\endxalignat
$$
Applying the last identity \mythetag{1.3} to \mythetag{1.4} and
to other identities \mythetag{1.3}, we derive
$$
\xalignat 3
&\nabla\bar{\boldsymbol\gamma}=0,
&&\nabla\bbd=0,
&&\nabla\bar{\bold H}=0.
\ \qquad
\endxalignat
$$
The general relativity (the Einstein's theory of gravity) is a theory
with zero torsion. Exactly for this case we have the following theorem.
\mytheorem{1.2} There is a unique metric connection $(\Gamma,\Alpha,
\bar{\Alpha})$ of the bundle of Dirac spinors $DM$ whose torsion 
$\bold T$ is zero.
\endproclaim
    The metric connection with zero torsion $\bold T=0$ is called 
the {\it Levi-Civita\/} connection. The proof of both 
theorems~\mythetheorem{1.1} and ~\mythetheorem{1.2} as well as some 
explicit formulas for the components of the Levi-Civita connection
can be found in \mycite{3}.\par
    A massive spin 1/2 particle is described by a wave-function which
is a smooth spinor field $\boldsymbol\psi$. In order to get a coordinate
representation of this field we choose two frames $(U,\,\boldsymbol
\Upsilon_0,\,\boldsymbol\Upsilon_1,\,\boldsymbol\Upsilon_2,\,\boldsymbol
\Upsilon_3)$ and $(U,\,\boldsymbol\Psi_1,\,\boldsymbol\Psi_2,\,\boldsymbol
\Psi_3,\,\boldsymbol\Psi_4)$ with common domain $U$. The first of these
two frames is given by four smooth vector fields $\boldsymbol \Upsilon_0$,
$\boldsymbol\Upsilon_1$, $\boldsymbol\Upsilon_2$, and $\boldsymbol
\Upsilon_3$ linearly independent at each point $p\in U$. The other frame
is formed by four spinor fields $\boldsymbol\Psi_1$, $\boldsymbol\Psi_2$,
$\boldsymbol\Psi_3$, $\boldsymbol\Psi_4$ also linearly independent at
each point $p\in U$. Having these two frames and taking their dual and
Hermitian conjugate frames, one easily get the coordinate representation
for an arbitrary spin-tensorial field. Using the components of the 
wave function $\boldsymbol\psi$ in the frame $(U,\,\boldsymbol\Psi_1,
\,\boldsymbol\Psi_2,\,\boldsymbol\Psi_3,\,\boldsymbol\Psi_4)$ we write 
the following action integral for this field:
$$
\gathered
S=i\,\hbar
\int\sum^4_{a=1}\sum^4_{\bar a=1}\sum^4_{b=1}\sum^3_{q=0}
D_{a\bar a}\,
\gamma^{\kern 0.5pt aq}_b
\frac{\overline{\psi^{\kern 0.5pt\bar a}}
\ \nabla_{\!q}\psi^{\kern 0.5pt\lower 1.2pt\hbox{$\ssize b$}}
-\psi^{\kern 0.5pt\lower 1.2pt\hbox{$\ssize b$}}\ \nabla_{\!q}
\overline{\psi^{\kern 0.5pt\bar a}}}{2}\,dV\,-\\
\vspace{2ex}
-\,m\,c\int\sum^4_{a=1}\sum^4_{\bar a=1}D_{a\bar a}\,
\overline{\psi^{\kern 0.5pt\bar a}}\,
\psi^{\kern 0.5pt\lower 1.2pt\hbox{$\ssize a$}}\,dV.
\endgathered
\mytag{1.5}
$$
Through $\hbar$ in \mythetag{1.5} we denote the Planck
constant\footnote{\ These data are taken from the NIST site
\myhref{http://physics.nist.gov/cuu/Constants}
{http:/\negskp/physics.nist.gov/cuu/Constants}.}, while $c$ 
is the speed of light:
\adjustfootnotemark{-1}
$$
\xalignat 2
&\hbar\approx 1.05457168\cdot 10^{-27}
\text{\it erg}\cdot\text{\it\! sec},
&&c\approx 2.99792458\cdot 10^{10}\,\text{\it cm}/\text{\it\! sec}.
\endxalignat
$$
The constant $m$ in \mythetag{1.5} is the mass of a particle. By $dV$ 
in \mythetag{1.5} we denote the $4$-dimensional volume element induced 
by the metric $\bold g$. In local coordinates $x^0,\,x^1\,x^2,\,x^3$ 
within the domain $U\subset M$ it is written as follows:
$$
\hskip -2em
dV=\sqrt{-\det\bold g\,}\,d^{\kern 0.5pt 4}\kern -0.5pt x=
\sqrt{-\det\bold g\,}\ dx^0\!\wedge dx^1\!\wedge dx^2\!\wedge dx^3.
\mytag{1.6}
$$
Though $D_{a\bar a}$, $\gamma^{\kern 0.5pt an}_b$, $\psi^{\kern 0.5pt
\lower 1.2pt\hbox{$\ssize a$}}$, and $\psi^{\kern 0.5pt\lower 1.2pt
\hbox{$\ssize b$}}$ are the components of complex fields, the integral
\mythetag{1.5} is a real quantity. This fact is proved with the use 
of the identities
$$
\xalignat 2
&\hskip -2em
\overline{D_{a\bar a}}=D_{\bar aa},
&&\sum^4_{a=1}D_{a\bar a}\,\gamma^{\kern 0.5pt aq}_b
=\sum^4_{\bar s=1}D_{\kern -0.5pt b\kern 0.5pt\bar s}\,
\overline{\gamma^{\kern 0.5pt\bar sq}_{\kern 0.5pt\bar a}}.
\qquad
\mytag{1.7}
\endxalignat 
$$
Then remember that the metric connection is a real connection. Therefore,
we have
$$
\hskip -2em
\nabla_{\!q}\overline{\psi^{\kern 0.5pt\bar a}}=
\overline{\nabla_{\!q}\psi^{\kern 0.5pt\bar a}}.
\mytag{1.8}
$$
Taking into account \mythetag{1.3}, \mythetag{1.4}, \mythetag{1.7},
and \mythetag{1.8}, one can easily derive $\overline{S}=S$ for the 
action integral \mythetag{1.5}.\par
     Applying the extremal action principle to the action integral
\mythetag{1.5}, we derive the following differential equation for the
components of the spinor wave-function $\boldsymbol\psi$\,:
$$
\hskip -2em
i\,\hbar\,\sum^4_{b=1}\sum^3_{q=0}
\gamma^{\kern 0.5pt aq}_b\,\nabla_{\!q}
\psi^{\kern 0.5pt\lower 1.2pt\hbox{$\ssize b$}}
-m\,c\ \psi^{\kern 0.5pt\lower 1.2pt\hbox{$\ssize a$}}=0.
\mytag{1.9}
$$
The equation \mythetag{1.9} is the {\it Dirac equation} for a spin 
1/2 particle with the rest mass $m$. Conservation laws for relativistic 
field equations are formulated in terms of currents. The vector-field
$\bold J$ with the components
$$
\hskip -2em
J^{\kern 0.5pt\lower 1.2pt\hbox{$\ssize q$}}
=c\,\sum^4_{a=1}\sum^4_{\bar a=1}\sum^4_{b=1}D_{a\bar a}
\,\gamma^{\kern 0.5pt aq}_b\,\overline{\psi^{\kern 0.5pt\bar a}}
\,\psi^{\kern 0.5pt\lower 1.2pt\hbox{$\ssize b$}}
\mytag{1.10}
$$
is a current for the Dirac equation \mythetag{1.9}. The conservation
law is written as 
$$
\hskip -2em
\divr\bold J=\sum^3_{q=0}\nabla_{\!q}J^{\kern 0.5pt\lower 1.2pt
\hbox{$\ssize q$}}=0.
\mytag{1.11}
$$
In the case of the current \mythetag{1.10} the conservation law 
\mythetag{1.11} is derived from the Dirac equation \mythetag{1.9}
with the use of the identities \mythetag{1.3} and \mythetag{1.4}.
\par
    The Dirac current \mythetag{1.10} is a real vector-field:
$\tau(\bold J)=\bold J$. Indeed, using the identities \mythetag{1.7}
and applying them to \mythetag{1.10}, one easily derives
$$
\overline{J^q}=J^{\kern 0.5pt\lower 1.2pt\hbox{$\ssize q$}}.
$$
Moreover, the vector-field $\bold J$ is composed by time-like vectors.
In order to prove this fact we calculate the following quantity:
$$
\hskip -2em
g(\bold J,\,\bold J)=
\sum^3_{p=0}\sum^3_{q=0}g_{p\kern 0.5pt q}
\,J^{\kern 0.5pt\lower 1.2pt\hbox{$\ssize p$}}
\,J^{\kern 0.5pt\lower 1.2pt\hbox{$\ssize q$}}.
\mytag{1.12}
$$
The quantity $g(\bold J,\,\bold J)$ in the left hand side of
\mythetag{1.12} is a scalar invariant of the vector $\bold J$.
Its value does not depend on a frame choice. Let's assume for
a while that $(U,\,\boldsymbol\Upsilon_0,\,\boldsymbol\Upsilon_1,
\,\boldsymbol\Upsilon_2,\,\boldsymbol\Upsilon_3)$ and $(U,\,\boldsymbol
\Psi_1,\,\boldsymbol\Psi_2,\,\boldsymbol\Psi_3,\,\boldsymbol\Psi_4)$
form a canonically associated frame pair such that $(U,\,\boldsymbol
\Upsilon_0,\,\boldsymbol\Upsilon_1,\,\boldsymbol\Upsilon_2,\,\boldsymbol
\Upsilon_3)$ is a positively polarized right orthonormal frame in $TM$
and $(U,\,\boldsymbol\Psi_1,\,\boldsymbol\Psi_2,\,\boldsymbol\Psi_3,
\,\boldsymbol\Psi_4)$ is a canonically orthonormal chiral frame in
$DM$ (see the diagram \thetag{5.12} in \mycite{3} for more details). 
In such a frame pair the Dirac $\gamma$-field is represented by the
following standard Dirac matrices:
$$
\xalignat 2
&\hskip -2em
\gamma^{a0}_b=\Vmatrix 0&0&1&0\\0&0&0&1\\1&0&0&0\\0&1&0&0\endVmatrix,
&&\gamma^{a1}_b=\Vmatrix 0&0&0&-1\\0&0&-1&0\\0&1&0&0\\1&0&0&0\endVmatrix,
\quad\\
\vspace{-1.5ex}
&&&\mytag{1.13}\\
\vspace{-1.5ex}
&\hskip -2em
\gamma^{a2}_b=\Vmatrix 0&0&0&i\\0&0&-i&0\\0&-i&0&0\\i&0&0&0\endVmatrix,
&&\gamma^{a3}_b=\Vmatrix 0&0&-1&0\\0&0&0&1\\1&0&0&0\\0&-1&0&0\endVmatrix.
\quad
\endxalignat
$$
Here $a$ stands for a raw number, while $b$ is a column number. For the
chirality operator and the Dirac form in such a frame pair we have
$$
\xalignat 2
&\hskip -2em
H^{\kern 0.5pti}_{\kern -0.5pt j}=\Vmatrix 1 & 0 & 0 & 0\\0 & 1 & 0 & 0\\
0 & 0 & -1 & 0\\0 & 0 & 0 & -1\endVmatrix,
&&D_{i\bar j}=\Vmatrix 0 & 0 & 1 & 0\\0 & 0 & 0 & 1\\
1 & 0 & 0 & 0\\0 & 1 & 0 & 0\endVmatrix.
\mytag{1.14}
\endxalignat
$$
For the metric tensor $\bold g$ and the spin-metric tensor $\bold d$ 
in such a frame pair we have
$$
\xalignat 2
&\hskip -2em
g_{ij}=\Vmatrix 1 & 0 & 0 & 0\\0 & -1 & 0 & 0\\
0 & 0 & -1 & 0\\0 & 0 & 0 & -1\endVmatrix,
&&d_{ij}=\Vmatrix 0 & 1 & 0 & 0\\-1 & 0 & 0 & 0\\
0 & 0 & 0 & -1\\0 & 0 & 1 & 0\endVmatrix.\quad
\mytag{1.15}
\endxalignat
$$
{\bf A remark}. Note that the matrix $\gamma^0$ in \mythetag{1.13}
coincides with the matrix $D$ in \mythetag{1.14}. For this reason
in many books $\gamma^0$ is used instead of $D$ (see \S~21 in
\mycite{4}, see section~14.1 in \mycite{5}, see section~5.4 in
\mycite{6}, see section~2.5 in \mycite{7}, and see section~6.3 
in \mycite{8}). This usage contradicts the spin-tensorial nature 
of the fields $\boldsymbol\gamma$ and $\bold D$ because the formulas
\mythetag{1.13}, \mythetag{1.14}, \mythetag{1.15} and the equality
$\gamma^0=D$ are highly frame-specific. I think the use of this equality
without indicating explicitly its restricted scope is misleading for many
generations of readers of the above very famous books.\par
    Returning back to the formula \mythetag{1.10} and applying the 
formulas \mythetag{1.13} and \mythetag{1.14} to it, we derive the 
following formulas for the current components:
$$
\align
\hskip -2em
J^{\kern 0.5pt\lower 1.2pt\hbox{$\ssize 0$}}
&=c^{\kern 0.5pt 2}\bigl(\,\psi^{\kern 0.5pt\lower
1.2pt\hbox{$\ssize 1$}}\,\overline{\psi^{\kern 0.5pt 1}}+\psi^{\kern
0.5pt\lower 1.2pt\hbox{$\ssize 2$}}\,\overline{\psi^{\kern 0.5pt 2}}
+\psi^{\kern 0.5pt\lower 1.2pt\hbox{$\ssize 3$}}\,\overline{\psi^{\kern
0.5pt 3}}+\psi^{\kern 0.5pt\lower 1.2pt\hbox{$\ssize 4$}}
\,\overline{\psi^{\kern 0.5pt 4}}\,\bigr),\qquad
\mytag{1.16}\\
\vspace{1ex}
\hskip -2em
J^{\kern 0.5pt\lower 1.2pt\hbox{$\ssize 1$}}
&=c^{\kern 0.5pt 2}\bigl(\,\psi^{\kern 0.5pt\lower 1.2pt\hbox{$\ssize 1$}}
\,\overline{\psi^{\kern 0.5pt 2}}+\psi^{\kern 0.5pt\lower 1.2pt
\hbox{$\ssize 2$}}\,\overline{\psi^{\kern 0.5pt 1}}-\psi^{\kern 0.5pt
\lower 1.2pt\hbox{$\ssize 3$}}\,\overline{\psi^{\kern 0.5pt 4}}
-\psi^{\kern 0.5pt\lower 1.2pt\hbox{$\ssize 4$}}
\,\overline{\psi^{\kern 0.5pt 3}}\,\bigr),\qquad\\
\hskip -2em
J^{\kern 0.5pt\lower 1.2pt\hbox{$\ssize 2$}}
&=i\,c^{\kern 0.5pt 2}\bigl(\,\psi^{\kern 0.5pt\lower 1.2pt\hbox{$\ssize 1$}}
\,\overline{\psi^{\kern 0.5pt 2}}-\psi^{\kern 0.5pt\lower 1.2pt
\hbox{$\ssize 2$}}\,\overline{\psi^{\kern 0.5pt 1}}-\psi^{\kern 0.5pt
\lower 1.2pt\hbox{$\ssize 3$}}\,\overline{\psi^{\kern 0.5pt 4}}
+\psi^{\kern 0.5pt\lower 1.2pt\hbox{$\ssize 4$}}
\,\overline{\psi^{\kern 0.5pt 3}}\,\bigr),\qquad
\mytag{1.17}\\
J^{\kern 0.5pt\lower 1.2pt\hbox{$\ssize 0$}}
&=c^{\kern 0.5pt 2}\bigl(\,\psi^{\kern 0.5pt\lower 1.2pt\hbox{$\ssize 1$}}
\,\overline{\psi^{\kern 0.5pt 1}}-\psi^{\kern 0.5pt\lower 1.2pt
\hbox{$\ssize 2$}}\,\overline{\psi^{\kern 0.5pt 2}}-\psi^{\kern 0.5pt
\lower 1.2pt\hbox{$\ssize 3$}}\,\overline{\psi^{\kern 0.5pt 3}}
+\psi^{\kern 0.5pt\lower 1.2pt\hbox{$\ssize 4$}}
\,\overline{\psi^{\kern 0.5pt 4}}\,\bigr).\qquad
\endalign
$$
Substituting \mythetag{1.16} and \mythetag{1.17} into
\mythetag{1.12}, we derive
$$
\hskip -2em
\gathered
g(\bold J,\,\bold J)=4\,c^{\kern 0.5pt 2}\bigl(\,\psi^{\kern 0.5pt
\lower 1.2pt\hbox{$\ssize 1$}}\,\overline{\psi^{\kern 0.5pt 1}}\,
\psi^{\kern 0.5pt\lower 1.2pt\hbox{$\ssize 3$}}\,\overline{\psi^{\kern
0.5pt 3}}+\psi^{\kern 0.5pt\lower 1.2pt\hbox{$\ssize 2$}}\,
\overline{\psi^{\kern 0.5pt 2}}\,\psi^{\kern 0.5pt\lower 1.2pt
\hbox{$\ssize 4$}}\,\overline{\psi^{\kern 0.5pt 4}}\,+\\
+\,\psi^{\kern 0.5pt\lower 1.2pt\hbox{$\ssize 1$}}\,\overline{\psi^{\kern
0.5pt 2}}\,\psi^{\kern 0.5pt\lower 1.2pt\hbox{$\ssize 4$}}
\,\overline{\psi^{\kern0.5pt 3}}+\psi^{\kern 0.5pt\lower 1.2pt
\hbox{$\ssize 2$}}\,\overline{\psi^{\kern 0.5pt 1}}\,\psi^{\kern 0.5pt
\lower 1.2pt\hbox{$\ssize 3$}}\,\overline{\psi^{\kern 0.5pt 4}}\,\bigr).
\endgathered
\mytag{1.18}
$$
Using the well-know inequality $2\,|x\,y|\leqslant |x|^2+|y|^2$, we get
the following inequality:
$$
\bigl|\psi^{\kern 0.5pt\lower 1.2pt\hbox{$\ssize 1$}}\,
\overline{\psi^{\kern 0.5pt 2}}\,\psi^{\kern 0.5pt\lower 
1.2pt\hbox{$\ssize 4$}}\,\overline{\psi^{\kern0.5pt 3}}
+\psi^{\kern 0.5pt\lower 1.2pt\hbox{$\ssize 2$}}\,
\overline{\psi^{\kern 0.5pt 1}}\,\psi^{\kern 0.5pt
\lower 1.2pt\hbox{$\ssize 3$}}\,\overline{\psi^{\kern 0.5pt 4}}
\bigr|\leqslant 2\,\bigl|\psi^{\kern 0.5pt\lower 1.2pt
\hbox{$\ssize 1$}}\,\psi^{\kern 0.5pt\lower 1.2pt\hbox{$\ssize 2$}}
\,\psi^{\kern 0.5pt\lower 1.2pt\hbox{$\ssize 3$}}\,\psi^{\kern 
0.5pt\lower 1.2pt\hbox{$\ssize 4$}}\bigr|\leqslant\bigl|
\psi^{\kern 0.5pt\lower 1.2pt\hbox{$\ssize 1$}}\,
\psi^{\kern 0.5pt\lower 1.2pt\hbox{$\ssize 3$}}\bigr|^2
+\bigl|\psi^{\kern 0.5pt\lower 1.2pt\hbox{$\ssize 2$}}\,
\psi^{\kern 0.5pt\lower 1.2pt\hbox{$\ssize 4$}}\bigr|^2.
$$
Applying this inequality to \mythetag{1.18}, we see that
$g(\bold J,\,\bold J)\geqslant 0$. This means that $\bold J$
is a time-like vector. Moreover, in \mythetag{1.16} we see 
that $J^{\kern 0.5pt\lower 1.2pt\hbox{$\ssize 0$}}\geqslant 0$. 
If we remember that for deriving \mythetag{1.16} and \mythetag{1.17} 
we took a positively polarized right orthonormal frame
$(U,\,\boldsymbol\Psi_1,\,\boldsymbol\Psi_2,\,\boldsymbol\Psi_3,
\,\boldsymbol\Psi_4)$ in $SM$, then from $J^{\kern 0.5pt\lower 1.2pt
\hbox{$\ssize 0$}}\geqslant 0$ we conclude that $\bold J$ is a time-like
vector from the interior of a Future light cone. Thus, we have proved the
following well-known theorem.
\mytheorem{1.3} The Dirac current \mythetag{1.10} for a massive
spin 1/2 particle is a time-like vector-field directed to the 
Future.
\endproclaim
\head
2. Normalization condition\\
for a single particle wave function.
\endhead
\parshape 21 0pt 360pt 0pt 360pt 0pt 360pt 0pt 360pt 0pt 360pt 0pt 360pt 
0pt 360pt 0pt 360pt 180pt 180pt 180pt 180pt 180pt 180pt 
180pt 180pt 180pt 180pt 180pt 180pt 180pt 180pt 
180pt 180pt 180pt 180pt 180pt 180pt 180pt 180pt 180pt 180pt 
0pt 360pt
    Let's remember that wave functions in quantum mechanics are usually
normalized. In the non-relativistic theory scalar wave functions of bound
states are normalized to unity by means of the following integral
(see \mycite{9}):
$$
\hskip -2em
\int|\psi|^2\,d^{\kern 0.5pt 3}\kern -0.5pt x=1.
\mytag{2.1}
$$
In the case of a Dirac particle in a non-flat space-time $M$ the integral
\mythetag{2.1} is senseless since there is no predefined $3$-dimensional
submanifold in $M$. \vadjust{\vskip 5pt\hbox to 0pt{\kern 10pt
\includegraphics{DiracEq01.eps}\hss}\vskip -5pt}However, using 
the Dirac current \mythetag{1.10} one can give a new sense to integrals
like \mythetag{2.1}. Let $S$ be some arbitrary space-like hypersurface 
in $M$. A part of such a hypersurface $S$ enclosed into a space-time 
cube is shown on Fig\.~2.1. At each point of $S$ there is a unique unit
normal vector $\bold n$ directed to the Future. Then the equality
$$
\hskip -2em
\int\limits_S g(\bold J,\bold n)\,dS=1
\mytag{2.2}
$$
is a proper normalization condition for the wave-function of a single
massive spin 1/2 particle. By $dS$ in \mythetag{2.2} we denote the
$3$-dimensional area element determined by the metric induced from $M$
to $S$. In local coordinates $u^1,\,u^2,\,u^3$ of $S$ it is given by
a formula similar to \mythetag{1.6}:
$$
dS=\sqrt{-\det\bold g\,}\,d^{\kern 0.5pt 3}\kern -0.5pt u=
\sqrt{-\det\bold g\,}\ du^1\!\wedge du^2\!\wedge du^3.
$$
The quantity $g(\bold J,\bold n)$ integrated in \mythetag{2.2} is the
scalar product of the Dirac current and the unit normal vector of $S$
calculated in the Minkowski metric $\bold g$:
$$
\hskip -2em
g(\bold J,\bold n)=\sum^3_{i=0}\sum^3_{j=0}g_{ij}
\,J^{\kern 0.5pt\lower 1.2pt\hbox{$\ssize i$}}
\,n^{\kern 0.5pt\lower 1.2pt\hbox{$\ssize j$}}.
\mytag{2.3}
$$
Since both $\bold J$ and $\bold n$ are time-like vectors directed to the
Future, the scalar product \mythetag{2.3} is a positive quantity.\par
     Assume that $\psi^{\kern 0.5pt\lower 1.2pt\hbox{$\ssize\bar a$}}$
and $\psi^{\kern 0.5pt\lower 1.2pt\hbox{$\ssize b$}}$ in \mythetag{1.10}
are the components of the wave function satisfying the Dirac equation
\mythetag{1.9}. Then the components of the current $\bold J$ satisfy
the differential equation \mythetag{1.11}. Assume that $S'$ is some
other space-like hypersurface in $M$ such that $S-S'$ is the boundary 
for some domain $\Omega$ (see Fig\.~2.1):
$$
\hskip -2em
S-S'=\partial\Omega.
\mytag{2.4}
$$
In this case from \mythetag{1.11} and \mythetag{2.4} we derive
$$
\hskip -2em
\int\limits_S g(\bold J,\bold n)\,dS-\int\limits_{S'} 
g(\bold J,\bold n)\,dS=
\int\limits_\Omega\divr\bold J\ dV=0.
\mytag{2.5}
$$ 
The equality \mythetag{2.5} means that the choice of the hypersurface 
$S$ in the normalization condition \mythetag{2.2} is inessential. This
normalization condition is preserved in time dynamics given by the Dirac
equation \mythetag{1.9}.\par
     Let $\boldsymbol\varphi$ and $\boldsymbol\psi$ be two different
wave functions corresponding to different quantum states of a massive
spin 1/2 particle. They both satisfy the Dirac equation \mythetag{1.9}.
By analogy to \mythetag{1.10} we define the current $\bold J(\boldsymbol
\varphi,\boldsymbol\psi)$ with the following components:
$$
\hskip -2em
J^{\kern 0.5pt\lower 1.2pt\hbox{$\ssize q$}}(\boldsymbol\varphi,
\boldsymbol\psi)=c\,\sum^4_{a=1}\sum^4_{\bar a=1}\sum^4_{b=1}D_{a\bar a}
\,\gamma^{\kern 0.5pt aq}_b\,\overline{\varphi^{\kern 0.5pt\bar a}}
\,\psi^{\kern 0.5pt\lower 1.2pt\hbox{$\ssize b$}}.
\mytag{2.6}
$$
Like the initial Dirac current \mythetag{1.10}, this current \mythetag{2.6}
satisfies the differential equation \mythetag{1.11}. Therefore, relying on
\mythetag{2.4} and \mythetag{2.5}, we define the pairing
$$
\hskip -2em
\bigl<\boldsymbol\varphi\,|\,\boldsymbol\psi\bigr>
=\int\limits_S g(\bold J(\boldsymbol\varphi,\boldsymbol\psi),
\bold n)\,dS.
\mytag{2.7}
$$
The pairing \mythetag{2.7} is a Hermitian pairing in the sense of the
following equality:
$$
\hskip -2em
\bigl<\boldsymbol\varphi\,|\,\boldsymbol\psi\bigr>=
\overline{\bigl<\boldsymbol\psi\,|\,\boldsymbol\varphi\bigr>}.
\mytag{2.8}
$$
Due to the theorem~\mythetheorem{1.3} it is a positive pairing:
$$
\hskip -2em
\bigl<\boldsymbol\psi\,|\,\boldsymbol\psi\bigr>=\Vert\boldsymbol\psi
\Vert^2\geqslant 0.
\mytag{2.9}
$$
Moreover, if $\Vert\boldsymbol\psi\Vert=0$, then $\boldsymbol
\psi=0$ almost everywhere on the hypersurface $S$ in the sense 
of the $3$-dimensional Lebesgue \pagebreak measure on $S$.\par
     The pairing \mythetag{2.7} is preserved in time dynamics 
determined by the Dirac equation \mythetag{1.9}. Due to \mythetag{2.8} 
and \mythetag{2.9} it defines the Hilbert space of quantum states of a
single massive spin 1/2 particle. We denote it $\Cal H_1$.
\head
3. Multiparticle wave-functions.
\endhead
    Let $\boldsymbol\psi_{[0]}(p),\,\boldsymbol\psi_{[1]}(p),\,
\boldsymbol\psi_{[2]}(p),\,\boldsymbol\psi_{[3]}(p),\,\ldots$, 
where $p\in M$, be a series of single particle wave-functions 
forming an orthonormal basis in the Hilbert space $\Cal H_1$:
$$
\hskip -2em
\bigl<\boldsymbol\psi_{[i]}\,|\,\boldsymbol\psi_{[j]}\bigr>=
\delta_{[ij]}.
\mytag{3.1}
$$
Note that in \mythetag{3.1} the indices are enclosed into the square 
brackets. This is done in order to distinguish them from tensorial and
spin-tensorial indices enumerating the components of wave-functions.
Let $p_{[1]},\,\ldots,\,p_{[n]}$ be $n$ points of the space-time
symbolizing the positions of $n$ particles. Then a multiparticle
wave-function can be constructed as a product of single particle 
wave-functions:
$$
\hskip -2em
\boldsymbol\psi_{[i_1]}(p_{[1]})\otimes\ldots\otimes
\boldsymbol\psi_{[i_n]}(p_{[n]}).
\mytag{3.2}
$$
In a coordinate form, i\.\,e\. upon choosing some frame pair $(U,\,
\boldsymbol\Upsilon_0,\,\boldsymbol\Upsilon_1,\,\boldsymbol
\Upsilon_2,\,\boldsymbol\Upsilon_3)$ and $(U,\,\boldsymbol\Psi_1,\,
\boldsymbol\Psi_2,\,\boldsymbol\Psi_3,\,\boldsymbol\Psi_4)$, the
wave-function \mythetag{3.2} is represented as
$$
\psi^{\,b_1}_{[i_1]}(p_{[1]})\cdot\ldots\cdot
\psi^{\,b_n}_{[i_n]}(p_{[n]}).
$$
Though the wave-function \mythetag{3.2} is a tensor product of 
$n$ spin-tensorial fields of the type $(1,0|0,0|0,0)$, it is not 
a spin-tensorial field itself since the multiplicands are 
spin-tensors at $n$ different points $p_{[1]},\,\ldots,\,p_{[n]}$.
The wave-function \mythetag{3.2} satisfies the Dirac equation
\mythetag{1.9} with respect to each its argument $p_s$:
$$
\sum^4_{b_s=1}\!\left(\!i\,\hbar\,\shave{\sum^3_{q=0}}
\gamma^{\kern 0.5pt aq}_{b_s}(p_s)\,\nabla_{\!q}^{[s]}
-m\,c\ \delta^a_{b_s}\!\right)\!
\psi^{\,b_1}_{[i_1]}(p_{[1]})\cdot\ldots\cdot
\psi^{\,b_s}_{[i_s]}(p_{[s]})\cdot\ldots\cdot
\psi^{\,b_n}_{[i_n]}(p_{[n]})=0.
$$
Lets consider some other multiparticle wave-function of the form
\mythetag{3.2}:
$$
\hskip -2em
\boldsymbol\psi_{[j_1]}(p_{[1]})\otimes\ldots\otimes
\boldsymbol\psi_{[j_n]}(p_{[n]}).
\mytag{3.3}
$$
Using the wave-functions \mythetag{3.2} and \mythetag{3.3} and 
applying the formula \mythetag{2.6} to them, we can define the 
following multicurrent:
$$
\hskip -2em
\bold J_{[i_1j_1]}(p_{[1]})\otimes\ldots\otimes
\bold J_{[i_nj_n]}(p_{[n]}).
\mytag{3.4}
$$
Here $\bold J_{[i_sj_s]}=\bold J(\boldsymbol\psi{[i_s]},
\boldsymbol\psi{[j_s]})$ for $s=1,\,\ldots,\,n$. One can integrate the
multicurrent \mythetag{3.4} over the Cartesian product of $n$ copies of
the hypersurface $S$ thus defining a pairing for multiparticle
wave functions:
$$
\hskip -2em
\bigl<\boldsymbol\psi_{[i_1]}\otimes\ldots\otimes\boldsymbol
\psi_{[i_n]}\,|\,\boldsymbol\psi_{[j_1]}\otimes\ldots\otimes
\boldsymbol\psi_{[j_n]}\bigr>=\prod^n_{s=1}\delta_{[i_sj_s]}.
\mytag{3.5}
$$
The formula \mythetag{3.5} shows that the wave-functions of the form
\mythetag{3.2} constitute an orthonormal basis in a Hilbert space
defined by the multicurrent \mythetag{3.4}. This Hilbert space is
denoted $\Cal H_n$. It is the tensor product of $n$ copies of $\Cal H_1$:
$$
\hskip -2em
\Cal H_n=\underbrace{\Cal H_1\otimes\ldots\otimes
\Cal H_1}_{\text{$n$ times}}.
\mytag{3.6}
$$\par
    Note that the Hilbert space \mythetag{3.6} is not a space of 
wave-functions for actual quantum states of $n$ particles. Wave-functions
of actual states should be symmetric for bosons and skew-symmetric for
fermions. In our case of spin 1/2 particles they are fermions. For this
reason we construct an actual wave-function as follows:
$$
\hskip -2em
\boldsymbol\psi_{[i_1\ldots\,i_n]}
=\sum_{\sigma\in\goth S_n}\!\!\frac{(-1)^\sigma}{\sqrt{n!}}
\ \boldsymbol\psi_{[i_{\sigma 1}]}(p_{[1]})\otimes\ldots\otimes
\boldsymbol\psi_{[i_{\sigma n}]}(p_{[n]}).
\mytag{3.7}
$$
In a coordinate form, i\.\,e\. upon choosing some frame pair $(U,\,
\boldsymbol\Upsilon_0,\,\boldsymbol\Upsilon_1,\,\boldsymbol
\Upsilon_2,\,\boldsymbol\Upsilon_3)$ and $(U,\,\boldsymbol\Psi_1,\,
\boldsymbol\Psi_2,\,\boldsymbol\Psi_3,\,\boldsymbol\Psi_4)$, the
wave-function \mythetag{3.7} is represented as
$$
\hskip -2em
\psi^{\,b_1\ldots\,b_n}_{[i_1\ldots\,i_n]}
=\sum_{\sigma\in\goth S_n}\!\!\frac{(-1)^\sigma}{\sqrt{n!}}
\ \psi^{\,b_1}_{[i_{\sigma 1}]}(p_{[1]})\otimes\ldots\otimes
\psi^{\,b_n}_{[i_{\sigma n}]}(p_{[n]}).
\mytag{3.8}
$$
By $\sigma$ in \mythetag{3.7} and \mythetag{3.8} we denote a 
transposition from the $n$-th symmetric group $\goth S_n$. Like
the wave-function \mythetag{3.7}, the wave-function \mythetag{3.7} 
satisfies the Dirac equation \mythetag{1.9} with respect to 
each its argument $p_s$:
$$
\sum^4_{b_s=1}\!\left(\!i\,\hbar\,\shave{\sum^3_{q=0}}
\gamma^{\kern 0.5pt aq}_{b_s}(p_s)\,\nabla_{\!q}^{[s]}
-m\,c\ \delta^a_{b_s}\!\right)\!
\psi^{\,b_1\ldots\,b_n}_{[i_1\ldots\,i_n]}(p_{[1]},\ldots,
p_{[n]})=0.\quad
\mytag{3.9}
$$
Due to \mythetag{3.9} one can apply the pairing \mythetag{3.5}
defined by means of the multicurrent \mythetag{3.4} to functions
of the form \mythetag{3.7}. As a result we get
$$
\hskip -2em
\bigl<\boldsymbol\psi_{[i_1\ldots\,i_n]}
\,|\,\boldsymbol\psi_{[j_1\ldots\,j_n]}
\bigr>=\prod^n_{s=1}\delta_{[i_sj_s]}.
\mytag{3.10}
$$
The formula \mythetag{3.10} means that the wave-functions of the form
\mythetag{3.7} constitute an orthonormal basis in a subspace of the 
Hilbert space \mythetag{3.6}. We denote this subspace through 
$\Cal H^{\text{\kern 0.5pt skew}}_n$. Note that $\Cal H^{\text{\kern 0.5pt
skew}}_1=\Cal H_1$. Let's denote $\Cal H_0=\Bbb C$ and $\boldsymbol
\Phi_0=1$. Then we consider the following direct sum of Hilbert spaces:
$$
\Cal H=\Cal H_0\oplus\bigoplus^\infty_{n=1}\Cal H^{\text{\kern 0.5pt 
skew}}_n.
\mytag{3.11}
$$
The pairing \mythetag{3.10} can be extended to \mythetag{3.11} so that
$$
\align
&\hskip -2em
\Cal H^{\text{\kern 0.5pt skew}}_n\!\perp \Cal H^{\text{\kern 0.5pt
skew}}_m\text{\ \ for all \ }n\neq m,\\
&\hskip -2em
\Cal H^{\text{\kern 0.5pt skew}}_n\!\perp \Cal H_0\text{\ \ for all 
\ }n,
\mytag{3.12}
\\
&\hskip -2em
|\Phi_0|=1.
\endalign
$$
Then due to \mythetag{3.12} the direct sum $\Cal H$ gains the structure 
of a Hilbert space. This is the Hilbert space of all multiparticle quantum
states of massive spin 1/2 particles described by the Dirac equation 
\mythetag{1.9}. The vector $\Phi_0\in\Cal H$ is called the {\it vacuum
vector\/} in the secondary quantization scheme (see \mycite{10}). 
According to this scheme the creation operators are introduced as 
follows:
$$
\align
&\hskip -2em
a^{\sssize +}_{[i_s]}\,\Phi_0=\boldsymbol\psi_{[i_s]},
\mytag{3.13}\\
&\hskip -2em
a^{\sssize +}_{[i_s]}\,\boldsymbol\psi_{[i_s]}=0,
\mytag{3.14}\\
&\hskip -2em
a^{\sssize +}_{[i_s]}\,\boldsymbol\psi_{[i_1\ldots\,i_{s-1}
i_{s+1}\ldots\,i_n]}=\boldsymbol\psi_{[i_1\ldots\,i_{s-1}i_s
i_{s+1}\ldots\,i_n]},
\mytag{3.15}\\
&\hskip -2em
a^{\sssize +}_{[i_s]}\,\boldsymbol\psi_{[i_1\ldots\,i_{s-1}i_s
i_{s+1}\ldots\,i_n]}=0,
\mytag{3.16}
\endalign
$$
where $i_1<\ldots<i_s<\ldots<i_n$. The annihilation operators are defined
as Hermitian conjugates for creation operators with respect to the pairing
\mythetag{3.10} extended to the Hilbert space \mythetag{3.11}. From
\mythetag{3.13}, \mythetag{3.14}, \mythetag{3.15}, and \mythetag{3.16}
one easily derives
$$
\align
&\hskip -2em
a_{[i_s]}\,\boldsymbol\psi_{[i_s]}=\Phi_0,
\mytag{3.17}\\
&\hskip -2em
a_{[i_s]}\,\Phi_0=0,
\mytag{3.18}\\
&\hskip -2em
a_{[i_s]}\,\boldsymbol\psi_{[i_1\ldots\,i_{s-1}i_s
i_{s+1}\ldots\,i_n]}=\boldsymbol\psi_{[i_1\ldots\,i_{s-1}
i_{s+1}\ldots\,i_n]},
\mytag{3.19}\\
&\hskip -2em
a_{[i_s]}\,\boldsymbol\psi_{[i_1\ldots\,i_{s-1}
i_{s+1}\ldots\,i_n]}=0.
\mytag{3.20}
\endalign
$$
Note that the creation and annihilation operators introduced by 
the formulas \mythetag{3.13}, \mythetag{3.14}, \mythetag{3.15},
\mythetag{3.16}, \mythetag{3.17}, \mythetag{3.18}, \mythetag{3.19},
and \mythetag{3.20} are constant operators in a constant Hilbert 
space $\Cal H$. This fact is not surprising. In the absence of
interaction terms in the action integral \mythetag{1.5} the quantum
states of Dirac particles remain unchanged in dynamics regardless 
to the number of particles we have.
\head
4. Some conclusions.
\endhead
     All results of this paper are known. In the case of the flat
Minkowski space they are broadly known. The main goal of the present
paper is not to claim a new result, but to emphasize the existence
of five basic spin-tensorial fields \mythetag{1.1} in the theory of
Dirac particles and to fix the novel notations for these basic fields 
(see the remark below the formula \mythetag{1.15}).
\head
5. Acknowledgments.
\endhead
    I am grateful to E.~G.~Neufeld who gave me the book \mycite{7}. I am
also grateful to V.~R.~Kudashev who gave me the book \mycite{6} many years
ago in exchange for a VINITI\,\footnotemark\ book with I.~M.~Krichever's
article.\footnotetext{\ VINITI is the Russian Institute for Scientific and
Technological Information ({\eightcyr Vserossi\ae\-ski{\ae} Institut
Nauchno{\ae} i Tehnichesko{\ae} Informatsii}). It publishes serial books
with selected and invited articles in various areas of science and technology.}
\par

\enddocument
\end